\documentclass{elsart3-1}

\usepackage{amssymb}
\usepackage{mathrsfs}
\usepackage[english,francais]{babel}
\usepackage{mathtools}

\newtheorem{theorem}{Theorem}[section]

\newtheorem{e-proposition}[theorem]{Proposition}
\newtheorem{corollary}[theorem]{Corollary}
\newtheorem{e-definition}[theorem]{Definition\rm}
\newtheorem{remark}{\it Remark\/}

\newtheorem{theoreme}{Th\'eor\`eme}[section]

\newtheorem{proposition}[theoreme]{Proposition}

\renewcommand{\theequation}{\arabic{equation}}
\def\og{\leavevmode\raise.3ex\hbox{$\scriptscriptstyle\langle\!\langle$~}}
\def\fg{\leavevmode\raise.3ex\hbox{~$\!\scriptscriptstyle\,\rangle\!\rangle$}}
\newcommand{\bqn}{\begin{eqnarray}}
\newcommand{\eqn}{\end{eqnarray}}
\newcommand{\bqns}{\begin{eqnarray*}}
\newcommand{\eqns}{\end{eqnarray*}}
\newcommand{\beq}{\begin{equation}}
\newcommand{\beqa} {\begin{array}{rl}}
\newcommand{\eeq}{\end{equation}}
\newcommand{\eeqa}{\end{array}}
\newcommand{\R}{\mathbb{R}}
\newcommand{\C}{\mathbb{C}}
\newcommand{\N}{\mathbb{N}}

\begin{document}
\begin{frontmatter}


\selectlanguage{english}
\title{A short remark on a Growth-Fragmentation Equation}


\selectlanguage{english}
\author[authorlabel]{M. Escobedo}
\ead{miguel.escobedo@ehu.es}

\address[authorlabel]{Departamento de Matem\'aticas, Universidad del Pa{\'\i}s Vasco (UPV/EHU). E--48080 Bilbao (SPAIN)}

\medskip
\begin{center}

\end{center}

\begin{abstract}
\selectlanguage{english}
An explicit solution for  a growth fragmentation equation with constant dislocation measure is obtained. In this example  the necessary  condition  for the  general results  in \cite{BW} about the existence of global solutions  in the so called self similar case is not satisfied. The solution is local  and blows up in finite time.
\end{abstract}
\end{frontmatter}
\selectlanguage{english}

\section{Introduction}
\label{I}
The purpose of this note is to present an explicit solution that blows up in finite time to the growth fragmentation equation
\bqn
\label{eq:croisfrag2}
&&\frac {\partial  u} {\partial t}(t,x) + \frac {\partial } {\partial x}\big(x^{1+\gamma } u (t,x)\big) +  x^\gamma u(t,x)=\int\limits_x^\infty \frac {1} {y} k_0\left(\frac {x} {y} \right)y^\gamma   u(t,y) dy,\,\,\,t>0, x>0\\
\label{S1EK0}
&&\gamma >0,\quad k_0(x)=\theta H(1-x),\,\,\,\, \theta>1,\,H: \hbox{Heaviside's function,}
\eqn 

Motivated by the study of growth-fragmentation stochastic processes \cite{BCK}, this  type of equation
was  considered  recently by J. Bertoin and A. R. Watson in  \cite{BW},   with the initial data 
\bqn
u(0, x)=\delta (x-1), \label{eq:croisfrag2data}
\eqn 
for  $\gamma\in \R$ and  $k_0$ a measure's density, with support contained in $[0, 1] $, that satisfies:
\bqn
\label{S0E1}
k_0(x)dx=k_0(1-x)dx,\,\,\,\forall x\in [1/2, 1); \quad\int  _{ [1/2, 1) }(1-x)^2k_0(x)dx=1.
\eqn

These equations have  proved to be  interesting for mathematical reasons (cf. \cite{BW}, \cite{DE}) and also because of the great variety of their applications in mathematical modeling (cf. \cite{BL}, \cite{DG}).

For $\gamma =0$, existence and uniqueness of non negative  solution to (\ref{eq:croisfrag2}),(\ref{eq:croisfrag2data}) is proved in   \cite{BW} under conditions (\ref{S0E1}) only. When $\gamma \not =0$ the  existence of a global solution  in  \cite{BW}   is proved with   the supplementary hypothesis
\bqn
&&\inf _{ s\ge 0 }\Phi (s)<0, \,\,\,\,\hbox{where:} \label{S1prop1B}\\
&&\Phi (s)= \left(K(s)+s-2\right),\,\,\,K(s)=\int _0^\infty x^{s-1}k_0(s)ds.\label{S1prop2}
\eqn

After the results in \cite{BW} and \cite{BS}, the importance of condition (\ref{S1prop1B}),(\ref{S1prop2}) is well established for  growth fragmentation processes, but it remains  to be better understood for the growth fragmentation equation.

We are considering in this note the simplest possible choice for $k_0$, given in  (\ref{S1EK0}). It  satisfies the condition (\ref{S0E1}), and is such that:
\bqn
&&K(s)=\frac {\theta} {s}\,\,\,\hbox{and}\,\,\,\,\Phi (s)=\frac {\theta} {s}+s-2\equiv \frac{(s-\sigma _1)(s-\sigma _2)}{s},\,\,\forall s\in \C;\,\,\Re e(s)>0, \label{S1EKS0}\\
&&\sigma _1=1-\sqrt{1-\theta},\,\,\,\sigma _2=1+\sqrt{1-\theta}.
\eqn
For $\theta\in (0, 1)$ the two roots of $\Phi (s)$ are positive real numbers and then condition  (\ref{S1prop1B}) is satisfied. 
But, for $\theta>1$, $\sigma _1$ and $\sigma _2$ are complex conjugated, then  $\inf _{ s>0 } \Phi (s)=2(\sqrt \theta -1)$ and  (\ref{S1prop1B}) is not satisfied.

Our main  result is the following Theorem, where $\mathscr{D}_1'$ denotes the set of distributions of order one.

\begin{theorem}
\label{S1MainTh1} 
For all $\gamma >0$ the measure on $\left(0, \gamma ^{-1}\right)\times (0, \infty)$, defined by:
\bqn
u(t, x)&=&u^S(t, x)+u^R(t, x) \label{S1Esol0}\\
u^S(t, x)&=&\left(1-\gamma  t\right)^{\frac {1} {\gamma }}\delta \left(x-\left(1-\gamma  t\right)^{-\frac {1} {\gamma }} \right)\label{S1Esol00}\\
u^R(t, x)&=&\theta\left(1-\gamma  t\right)^{\frac {2} {\gamma }}tF\left(1+\frac {\sigma _1} {\gamma }, 1+\frac {\sigma _2} {\gamma }, 2, \gamma  t\left(1+\left(\gamma  t-1\right)x^\gamma  \right) \right)H\left(1-\left(1-\gamma  t\right)^{\frac {1} {\gamma }}x\right),
\label{S1Esol}
\eqn
is non negative and satisfies,  the equation(\ref{eq:croisfrag2}),(\ref{S1EK0}) in $ \mathscr{D}_1'\left(\left(0, \gamma ^{-1}\right)\times (0, \infty)\right)$. It also satisfies  $u(t)\rightharpoonup \delta (x-1)$ in the weak sense of measures as  $t\to 0$.  
\end{theorem}

As a Corollary we deduce the following
\begin{corollary}
\label{S1Cor}
The  solution $u$ defined in (\ref{S1Esol0}),(\ref{S1Esol}) satisfies:
\bqn
\label{S1Ebw0}
\lim _{ \gamma  t\to 1^- }u(t, x)=\frac {\gamma \Gamma \left(\frac {2} {\gamma } \right)}
{  \Gamma \left(\frac {\sigma _1} {\gamma }  \right) \Gamma \left(\frac {\sigma _2} {\gamma }  \right)}(1+x^\gamma )^{-\frac {2} {\gamma }},\,\,\,\forall x>0
\eqn
and it blows up in finite time in the following sense:
\bqn
\forall r>1:\,\,\, \lim _{ t\to \gamma ^{-1} }(1-\gamma t)^{\frac {r-1} {\gamma }}\int _0^\infty x^r u(t, x)dx&=&\frac {\Gamma \left(\frac {r+1} {\gamma }\right)\Gamma \left(\frac {r-1} {\gamma }\right)}
{\Gamma \left(\frac {r+1-\sigma _1} {\gamma }\right)\Gamma \left(\frac {r+1-\sigma _2} {\gamma }\right)},\label{S1Ebw1} \\
 \lim _{ t\to \gamma ^{-1} }\frac {-1} {\log (1-\gamma t)}\int _0^\infty x u(t, x)dx&=&\frac {\Gamma \left(\frac {r} {\gamma }\right)}
{\Gamma \left(\frac {\sigma _1} {\gamma }\right)\Gamma \left(\frac {\sigma _2} {\gamma }\right)},\label{S1Ebw2}\\
\forall r\in (0, 1):\,\,\,\,
\lim _{ t\to \gamma ^{-1} }\int _0^\infty x^r u(t, x)dx&=&\frac {\Gamma \left(\frac {r+1} {\gamma }\right)\Gamma \left(\frac {1-r} {\gamma }\right)}
{\Gamma \left(\frac {\sigma _1} {\gamma }\right)\Gamma \left(\frac {\sigma _2} {\gamma }\right)}.\label{S1Ebw3}
\eqn
\end{corollary}

The question of the possible  extension  for $t>\gamma ^{-1}$ is beyond the scope of this note. More general dislocation measures like $
k_m(x)=(x^m+(1-x)^m)H(1-x)$ for $m=1, 2, 3, \cdots$ may also be considered (cf.  \cite{E}), although  the solutions are not always so explicit. 

\section{Mellin variables.}
\label{2}
If $u$ were a suitable solution of (\ref{eq:croisfrag2}),(\ref{eq:croisfrag2data}), applying the  Mellin transform to both sides  of (\ref{eq:croisfrag2}) and  (\ref{eq:croisfrag2data}),  we  would obtain for $\mathcal M_u$, the Mellin transform of $u$:
\bqn
\frac {\partial } {\partial t} \mathcal M_u(t, s)&=&( K(s)+s-2)\mathcal M_u(t, s+\gamma ) \label{S2eqmellin}\\
\mathcal M_u(0, s )&=&1. \label{S2initial} 
\eqn

Solutions to  (\ref{S2eqmellin}),(\ref{S2initial}) may be obtained by a general method, based on Wiener Hopf arguments (cf. \cite{E} for details).  For a description and applications of that method the  reader may consult \cite{BZ}. However in our case the problem  (\ref{S2eqmellin}),(\ref{S2initial}) has a particularly simple explicit  solution.

If  $F(a, b,  c, z)$  denotes the Gauss hypergeometric function $_2F_1(a, b, c, z)$ (see for exemple \cite{AS}),  it follows from  the identities 15.2.1 and 15.3.3 in \cite{AS},  that  the  function:
\bqn
\label{S2solmellin}
\Omega (t, s)=F\left(\frac {s-\sigma _1} {\gamma }, \frac {s-\sigma _2} {\gamma }, \frac {s} {\gamma }, \gamma t\right)\equiv (1-\gamma t)^{\frac {2-s} {\gamma }}F\left(\frac {\sigma _1} {\gamma }, \frac {\sigma _2} {\gamma  }, \frac {s} {\gamma }, \gamma t \right)
\eqn
satisfies (\ref{S2eqmellin}),(\ref{S2initial}) for $t$ and $s$ in the domain of analyticity of $F\left(\frac {s-\sigma _1} {\gamma }, \frac {s-\sigma _2} {\gamma }, \frac {s} {\gamma }, \gamma t\right)$ such that $t\not = \gamma ^{-1}$. 

Our purpose is to define  the fonction $u$ as  the inverse Mellin transform of  $\Omega $, to prove that the Mellin transform of $u$ is $\Omega $ and then to  prove that $u$ solves (\ref{eq:croisfrag2}),(\ref{eq:croisfrag2data}).  

\section{The inverse Mellin transform of $\Omega (t, s)$.}
\label{S4}
We  first  show the following Proposition, where   $\mathscr M(0, \infty)$ denotes the space of non negative locally bounded measures on $(0, \infty)$.

\begin{proposition}
\label{S3P10}
For all $t\in \left(0, \gamma ^{-1} \right)$ the function $\Omega (t, s)$ defined in (\ref{S2solmellin}) has an inverse Mellin transform that belongs to   $\mathscr M(0, \infty)$, that we denote $u (t, x)=\mathcal M ^{-1}(\Omega )$, and that satisfies
\bqn
\label{S3E124}
\mathcal M_{u} (t, s)=\Omega (t, s),\,\,\,\forall s\in \C; \,\, \Re e(s)>0, \,\, t\in \left(0, \gamma ^{-1}\right).
\eqn
\end{proposition}
\textbf{Proof.} For $0<t<\gamma ^{-1}$ fixed, the hypergeometric function $F\left(\frac {\sigma _1} {\gamma }, \frac {\sigma _2} {\gamma  }, \frac {s} {\gamma }, \gamma t \right)$ is analytic in the domain $D=\{s\in \C; \Re e(s)>0\}$ and by 15.7.1 in \cite{AS}
$$
\left|F\left(\frac {\sigma _1} {\gamma }, \frac {\sigma _2} {\gamma }, \frac {s} {\gamma }, \gamma t\right)-1-\frac {2 t } {s }\right|\le C(1+|s|)^{-2},\,\,\, \forall s\in D.
$$
for some constant $C=C(t, \sigma _1, \sigma _2, \gamma )$. Then, by Theorem 11.10.1 in \cite{ML}, for all $t\in \left(0, \gamma ^{-1}\right)$, the function $\Omega(t)$ has  an inverse Mellin transform $u(t)\in \mathscr M(0, \infty)$, given by
$$
u(t, x)=\frac {1} {2i\pi }\int  _{ \Re e (s)=s_0 }x^{-s}\Omega (t, s)ds,
$$
for an arbitrary $s_0>0$, and  such that for all $s\in D$, $\mathcal M _{ u(t) }(s)=\Omega (t, s)$. \qed

We may obtain now the explicit expression of $u$.  Let us prove first the following Proposition.
\begin{proposition}
\label{S4P154}
Suppose $\sigma _1\in \C$, $\sigma _2\in \C$,  $\gamma >0$ and $t\in \left(0, \gamma ^{-1} \right)$ and define the function
\bqn
v(t, x)=F\left(1+\frac {\sigma _1} {\gamma }, 1+\frac {\sigma _2} {\gamma }, 2, \gamma t\left(1+(\gamma t-1)x^\gamma  \right) \right)H\left(1-(1-\gamma t)^{\frac {1} {\gamma }}x\right)
\eqn
for $x>0$. Then, for all $t\in \left(0, \gamma ^{-1}\right)$, the Mellin transform of $v$ is:
\bqn
\label{S2Mellinv}
\mathcal M_v(t, s)= (1-\gamma t)^{-\frac {s} {\gamma }}\frac {F\left(\frac {\sigma _1} {\gamma }, \frac {\sigma _2} {\gamma }, \frac {s} {\gamma }, \gamma t\right)-1} 
{\theta t},\,\,\,\forall s\in \C;\,\,\,\Re e(s)>0.
\eqn
\end{proposition}

\textbf{Proof.} Since $\gamma >0$ and  $1-\gamma t>0$ it follows that for all $x>0$,  $\left(1+(\gamma t-1)x^\gamma  \right)<1$. Then:
\bqn
&&F\left(1+\frac {\sigma _1} {\gamma }, 1+\frac {\sigma _2} {\gamma }, 2, \gamma t\left(1+(\gamma t-1)x^\gamma  \right) \right)=
\sum_{ n=0 } ^\infty \frac {\Gamma \left(1+\frac {\sigma _1} {\gamma }+n \right)\Gamma \left( 1+\frac {\sigma _2} {\gamma }+n\right)\Gamma \left( 2\right)(\gamma t)^n} {\Gamma \left( 1+\frac {\sigma _1} {\gamma }\right)\Gamma \left( 1+\frac {\sigma _2} {\gamma }\right)\Gamma \left(2+n \right)\Gamma \left(n+1 \right)}\times \nonumber  \\
&&\hskip 11cm  \times(1+(\gamma t-1)x^\gamma )^n, \label{S4hyper}\\
&&\mathcal M_v(t, s)=\int_0^\infty  v(t, x)x^{s-1}dx=\sum_{ n=0 } ^\infty \frac {\Gamma \left(1+\frac {\sigma _1} {\gamma }+n \right)\Gamma \left( 1+\frac {\sigma _2} {\gamma }+n\right)\Gamma \left( 2\right)(\gamma t)^n} {\Gamma \left( 1+\frac {\sigma _1} {\gamma }\right)\Gamma \left( 1+\frac {\sigma _2} {\gamma }\right)\Gamma \left(2+n \right)\Gamma \left(n+1 \right)}\times \nonumber\\
&&\hskip 8cm \times \int_0^{(1-\gamma t)^{-\frac {1} {\gamma }}}(1+(\gamma t-1)x^\gamma )^n x^{s-1}dx.\nonumber
\eqn
A straightforward calculation gives, using $\gamma >0$:
$$
\int_0^{(1-\gamma t)^{-\frac {1} {\gamma }}}(1+(\gamma t-1)x^\gamma )^n x^{s-1}dx=(1-\gamma t)^{-\frac {s} {\gamma }}\frac {\Gamma (n+1)\Gamma \left( \frac {s} {\gamma }\right)} {\gamma \Gamma \left(1+\frac {s} {\gamma }+n \right)}.\nonumber
$$
and then,  for all $s\in \C$ such that  $-s/\gamma\not \in \N $: 
\bqns
\mathcal M_v(t, s)&=&(1-\gamma t)^{-\frac {s} {\gamma }}\sum_{ n=0 } ^\infty \frac {\Gamma \left(1+\frac {\sigma _1} {\gamma }+n \right)\Gamma \left( 1+\frac {\sigma _2} {\gamma }+n\right)\Gamma \left( 2\right)(\gamma t)^n} {\Gamma \left( 1+\frac {\sigma _1} {\gamma }\right)\Gamma \left( 1+\frac {\sigma _2} {\gamma }\right)\Gamma \left(2+n \right)\Gamma \left(n+1 \right)}\frac {\Gamma (n+1)\Gamma \left( \frac {s} {\gamma }\right)} {\gamma \Gamma \left(1+\frac {s} {\gamma }+n \right)}\\
&=& (1-\gamma t)^{-\frac {s} {\gamma }}\frac {F\left(\frac {\sigma _1} {\gamma }, \frac {\sigma _2} {\gamma }, \frac {s} {\gamma }, \gamma t\right)-1} {\sigma _1\sigma _2t}.
\eqns 
\qed

The next Corollary follows  from Proposition \ref{S4P154} and Theorem  11.10.1 in \cite{ML} on the uniqueness of the inverse Mellin transform:
\begin{corollary}
\label{coro1}
For all $\sigma _1\in \C$, $\sigma _2\in \C$,  suppose that $\gamma >0$, $0<\gamma t<1$ and let $\omega $ be the measure:
\bqns
\omega (t, x)=(1-\gamma t)^{\frac {1} {\gamma }}\delta \left(x-(1-\gamma t)^{-\frac {1} {\gamma }} \right)+ \sigma_1\sigma _2 t(1-\gamma t)^{\frac {2} {\gamma }}v(t, x).
\eqns
Then, for all $t\in \left(0, \gamma ^{-1}\right)$:
$$
\mathcal M_\omega(t, s)=\Omega (t, s)\,\,\, \hbox{for all}\,\,\,s\in \C,\, \Re e(s) s>0
$$
and $u(t)=\omega (t)$ for all $t\in \left(0, \gamma ^{-1}\right)$.
\end{corollary}

We may prove  now our  main result.

\noindent
\textbf{Proof of Theorem \ref{S1MainTh1}.} It is easy to  check that $u(t)\rightharpoonup \delta (x-1)$ as $t\to 0$ in the wek sense of measures.
We already  know that $\Omega (t, s) =\mathcal M_u (t, s)$ solves the problem (\ref{S2eqmellin}), (\ref{S2initial}).
Applying the inverse Mellin tranform to both sides of the equation  (\ref{S2eqmellin}) we deduce the following equation in 
 $\mathscr D'\left(\left(0, \gamma ^{-1}\right)\times (0, \infty)\right)$
\bqn
\label{S4E99}
\frac {\partial u} {\partial t} (t, x)=\frac {1} {2i\pi }\int  _{ \Re e s=\sigma _0 }\left( \frac {\theta} {s}+(s-1)-1\right)\mathcal M_u (t, s+\gamma )x^{-s}ds.
\eqn
We consider now each of the terms in the right and side separately. Since $\sigma _0>0$, $\gamma >0$, using that $\mathcal M_u (t, s )=\Omega (t, s)$ for all $\Re e(s)>0$ we write $\gamma >0$:
\bqn
\frac {1} {2i\pi }\int  _{ \Re e s=\sigma _0 }\mathcal M_u (t, s+\gamma )x^{-s}ds&=&\frac {1} {2i\pi }\int  _{ \Re e s=\sigma _0 }\int _0^\infty u(t, y)y^{s+\gamma -1}dyx^{-s}ds \nonumber\\
&=&\int _0^\infty u(t, y)y^{\gamma -1}\left(\frac {1} {2i\pi }\int  _{ \Re e s=\sigma _0 }\left(\frac{x}{y}\right)^{-s}ds\right)=x^\gamma u(t, x).\label{S4E100}
\eqn
The second term in the right hand side of (\ref{S4E99})  is given by the classical  formula
\bqn
\label{S4E101}
\frac {1} {2i\pi }\int  _{ \Re e s=\sigma _0 }(s-1)\mathcal M_u (t, s+\gamma )x^{-s}ds=\frac {\partial } {\partial x}\left(\frac {1} {2i\pi }\int  _{ \Re e s=\sigma _0 }\mathcal M_u (t, s+\gamma )x^{1-s}ds \right)
\eqn
In  the las term in  the right hand side of (\ref{S4E99}) we write as above:
$$
\frac {1} {2i\pi }\int  _{ \Re e s=\sigma _0 }\frac {\theta} {s}\mathcal M_u (t, s+\gamma )x^{-s}ds=
\int _0^\infty u(t, y)\left(\frac {1} {2i\pi }\int  _{ \Re e s=\sigma _0 }\frac {\theta} {s}y^{s+\gamma -1}x^{-s}ds\right)dy.
$$
Using that for $\sigma _0>0$:

\[
\frac {1} {2i\pi }\int  _{ \Re e s=\sigma _0 }\frac {1} {s}y^{s+\gamma -1}x^{-s}ds=
\begin{dcases}
0, & \hbox{if}\,\,y<x \\
 y^{\gamma -1}, &  \hbox{if}\,\, y>x
\end{dcases}
\]
we deduce
\bqn
\label{S4E102}
\frac {1} {2i\pi }\int  _{ \Re e s=\sigma _0 }\frac {\theta} {s}\mathcal M_u (t, s+\gamma )x^{-s}ds=
\theta \int _x^\infty u(t, y)y^{\gamma -1}dy.
\eqn
Since $u\in  C\left( \left[0, \gamma ^{-1}\right); \mathscr M(0, \infty)\right)\cap  C^1\left( \left(0, \gamma ^{-1}\right); \mathscr D_1'(0, \infty)\right)$ it  follows from (\ref{S4E99})--(\ref{S4E102}) that both sides of the equation  (\ref{eq:croisfrag2}) are equal  in 
$C\left( \left(0, \gamma ^{-1}\right); \mathscr D_1'(0, \infty)\right)$ and then, for all $\varphi \in C_c^1\left(\left(0, \frac {1} {\gamma }\right)\times (0, \infty)\right)$:
\bqn
\label{S4EW}
\left\langle u_t(t)+\frac {\partial } {\partial x}\left(x^{\gamma +1}u(t)\right)+x^\gamma u(t),\varphi (t)\right\rangle=\theta\left\langle\int _x^\infty y^{\gamma -1}u(t, y)dy, \varphi (t)\right\rangle
\eqn
where
 $\left\langle \cdot, \cdot\right\rangle$ is the duality bracket between  $ \mathscr{D}_1'\left(\left(0, \gamma ^{-1}\right)\times (0, \infty)\right)$ and  
 $C_c^1\left(\left(0, \gamma ^{-1}\right)\times (0, \infty)\right)$.
 
Since 
$\sigma_1=\overline{\sigma _2}$, $
\Gamma \left(1+\frac {\sigma _1} {\gamma }+n\right)=\overline{\Gamma \left(1+\frac {\sigma _2} {\gamma }+n\right)}$ for all  $n\in \N$, and the positivity of $u$ follows.\qed
\vskip 0.3cm 
\textbf{Proof of Corollary \ref{S1Cor}.}
Properties (\ref{S1Ebw1})--(\ref{S1Ebw3}) follow from the explicit expressions of $\Omega (t, s)$ and   $u^R$ in (\ref{S1Esol}), (\ref{S2solmellin}), and  formulas 15.4 (ii) in \cite{O}.
\qed
\vskip 0.3cm 
\begin{remark}
Due to the  particular form of the measure $u$, it is easy to check that (\ref{S4EW}) is satisfied for all $\varphi \in C_c^1\left(\left[0, \frac {1} {\gamma }\right)\times (0, \infty)\right)$. We also  deduce from (\ref{S4EW}) that $u^S$ and $u^R$ satisfy:
\bqn
\frac {\partial u^S} {\partial t}+\frac {\partial } {\partial x}\left(x^{\gamma +1}u^S \right)+x^\gamma u^S=0,\,\,\,\hbox{in}\,\,\,\mathscr D  _{ 1 }'(0, \infty)
\eqn
and for all $t>0$ and a.e. $x>0$:
\bqn
&&\frac {\partial u^R} {\partial t}(t, x)+\frac {\partial } {\partial x}\left(x^{\gamma +1}u^R \right)(t, x)+x^\gamma u^R(t, x)=\theta (1-\gamma t)^{-\frac {\gamma -1} {\gamma }}
H\left(1-\left(1-\gamma  t\right)^{\frac {1} {\gamma }}x\right)+\nonumber\\
&&\hskip 10cm +\theta\int _x^\infty u^R(t, y)y^{\gamma -1}dy.
\eqn
\end{remark}
\section*{Acknowledgements}
The author acknowledges support from  DGES Grant  MTM2014-52347-C2-1-R and Basque Government Grant IT641-13.

\end{document}